\documentclass[11pt,twoside]{amsart}
\usepackage{graphicx}
\usepackage{amssymb}
\usepackage{amsmath}
\usepackage{amsthm}
\usepackage{amsfonts}
\usepackage{amsmath, amsfonts, amssymb}
\usepackage{multirow}
\usepackage{rotating}
\usepackage{multirow}
\usepackage{rotating}
\usepackage{multirow}
\usepackage{rotating}
\newtheorem{theorem}{Theorem}[section]
\newtheorem{corollary}[theorem]{Corollary}

\newtheorem{definition}[theorem]{Definition}
\newtheorem{remark}[theorem]{\bf{Remark}}

\numberwithin{equation}{section}

\begin{document}
	
	\begin{center}\small{In the name of Allah, Most Gracious, Most Merciful.}\end{center}
	\vspace{0.5cm}
	
\title{Two-dimensional left (right) unital algebras over algebraically closed fields and $\mathbb{R}$}

\author{H.Ahmed$^1$, U.Bekbaev$^2$, I.Rakhimov$^3$}

\thanks{{\scriptsize
emails: $^1$houida\_m7@yahoo.com; $^2$bekbaev@iium.edu.my; $^3$rakhimov@upm.edu.my.}}
\maketitle
\begin{center}
{\scriptsize \address{$^{1}$Department of Math., Faculty of Science, UPM, Selangor, Malaysia $\&$ \\ Depart. of Math., Faculty of Science, Taiz University, Taiz, Yemen}\\
  \address{$^2$Department of Science in Engineering, Faculty of Engineering, IIUM, Malaysia}\\
\address{$^3$Institute for Mathematical Research (INSPEM), UPM, Serdang, Selangor, Malaysia}}
\end{center}
\begin{abstract} In this paper we describe all, up to isomorphism, left unital, right unital and unital algebra structures on two-dimensional vector space over any algebraically closed field and $\mathbb{R}$. We tabulate the algebras with the units.
\end{abstract}
\section{Introduction}

The principal building blocks of our descriptions are derived from \cite{A, B1} as the authors have presented complete lists of isomorphism classes of all two-dimensional algebras over algebraically closed fields and $\mathbb{R}$, providing the lists of canonical representatives of their structure constant's matrices. We will speak of a unity element and unital algebras rather than identity element and algebras with identity. The latest lists of all unital associative algebras
in dimension two, three, four, and five are available in \cite{Peirce}, \cite{Arezina}, \cite{Gabriel} and \cite{Mazolla}, respectively. The lists of all
associative algebras (both unital and non-unital) in dimension two and three are presented in \cite{Fialowski, Rakhimov1}. In this paper we describe the isomorphism classes of two-dimensional left(right) unital algebras over arbitrary algebraically closed field and $\mathbb{R}$. Our approach is totally different than that of \cite{ Arezina,  Fialowski, Gabriel, Mazolla, Peirce, Rakhimov1}.

We consider left(right) unital algebras over algebraically closed fields of characteristic not $2,3$, characteristic $2$, characteristic $3$ and over $\mathbb{R}$ separately according to classification results of \cite{A, B1}. To the best knowledge of authors the descriptions of left(right) unital two-dimensional algebras over
algebraically closed fields and $\mathbb{R}$ have not been given yet.
\section{ Preliminaries}

Let $\mathbb{F}$ be any field, $A\otimes B$ stand for the Kronecker product consisting of blocks $(a_{ij}B),$ where  $A=(a_{ij})$, $B$ are matrices over $\mathbb{F}$.

Let $(\mathbb{A},\cdot)$ be $m$-dimensional algebra over $\mathbb{F}$ and $e=(e^1,e^2,...,e^m)$ its basis. Then the bilinear map $\cdot$ is represented by a matrix $A=(A^k_{ij})\in M(m\times m^2;\mathbb{F})$ as follows \[\mathbf{u}\cdot \mathbf{v}=eA(u\otimes v),\] for $\mathbf{u}=eu,\mathbf{v}=ev,$
where $u = (u_1, u_2, ..., u_m)^T,$  $v = (v_1, v_2, ..., v_m)^T$ are column coordinate vectors of $\mathbf{u}$ and $\mathbf{v},$ respectively.
The matrix $A\in M(m\times m^2;\mathbb{F})$ defined above is called the matrix of structural constants (MSC) of $\mathbb{A}$ with respect to the basis $e$. Further we assume that a basis $e$ is fixed and we do not make a difference between the algebra
$\mathbb{A}$ and its MSC $A$ (see \cite{B}).

If $e'=(e'^1,e'^2,...,e'^m)$ is another basis of $\mathbb{A}$, $e'g=e$ with $g\in G=GL(m;\mathbb{F})$, and  $A'$ is MSC of $\mathbb{A}$ with respect to $e'$ then it is known that
\begin{equation}\label{3}A'=gA(g^{-1})^{\otimes 2}\end{equation} is valid. Thus, the isomorphism of algebras $\mathbb{A}$ and $\mathbb{B}$ over $\mathbb{F}$ can be given in terms of MSC as follows.

\begin{definition} Two $m$-dimensional algebras $\mathbb{A}$, $\mathbb{B}$ over $\mathbb{F}$, given by
	their matrices of structure constants $A$, $B$, are said to be isomorphic if $B=gA(g^{-1})^{\otimes 2}$ holds true for some $g\in GL(m;\mathbb{F})$.\end{definition}

\begin{definition}
An element $\mathbf{1_\textit{L}}$ ($\mathbf{1_\textit{R}}$) of an algebra $\mathbb{A}$ is called a left (respectively, right) unit if
$ \mathbf{1_\textit{L}}\cdot\mathbf{u}= \mathbf{u}$ (respectively, $\mathbf{u}\cdot \mathbf{1_\textit{R}}= \mathbf{u}$) for all $\mathbf{u}\in  \mathbb{A}.$ An algebra with the left(right) unit element is said to be left(right) unital algebra, respectively. \end{definition}
\begin{definition}
An element $\mathbf{1}\in \mathbb{A}$ is said to be an unit element  if $ \mathbf{1}\cdot\mathbf{u}=\mathbf{u}\cdot\mathbf{1}= \mathbf{u}$ for all $\mathbf{u}\in  \mathbb{A}.$ In this case the algebra $\mathbb{A}$ is said to be unital. \end{definition}

Further we consider only the case $m=2$ and for the simplicity we use \[A=\left(\begin{array}{cccc} \alpha_1 & \alpha_2 & \alpha_3 &\alpha_4\\ \beta_1
& \beta_2 & \beta_3 &\beta_4\end{array}\right)\] for MSC, where
$\alpha_1, \alpha_2, \alpha_3, \alpha_4, \beta_1, \beta_2, \beta_3, \beta_4$ stand for any elements of $\mathbb{F}$.\\

Due to \cite{A} we have the following classification theorems according to $Char(\mathbb{F})\neq 2,3,$ $Char(\mathbb{F})=2$ and $Char(\mathbb{F})=3$ cases, respectively.
\begin{theorem}\label{thm1} Over an algebraically closed field $\mathbb{F}$ $(Char(\mathbb{F})\neq 2$ and $3)$, any non-trivial $2$-dimensional algebra is isomorphic to only one of the following algebras listed by their matrices of structure constants:
	\begin{itemize}
	\item $A_{1}(\mathbf{c})=\left(
	\begin{array}{cccc}
	\alpha_1 & \alpha_2 &\alpha_2+1 & \alpha_4 \\
	\beta_1 & -\alpha_1 & -\alpha_1+1 & -\alpha_2
	\end{array}\right),\ \mbox{where}\ \mathbf{c}=(\alpha_1, \alpha_2, \alpha_4, \beta_1)\in \mathbb{F}^4,$
	\item $A_{2}(\mathbf{c})=\left(
	\begin{array}{cccc}
	\alpha_1 & 0 & 0 & 1 \\
	\beta _1& \beta _2& 1-\alpha_1&0
	\end{array}\right)\simeq \left(
	\begin{array}{cccc}
	\alpha_1 & 0 & 0 & 1 \\
	-\beta _1& \beta _2& 1-\alpha_1&0
	\end{array}\right),\ \mbox{where}\ \mathbf{c}=(\alpha_1, \beta_1, \beta_2)\in \mathbb{F}^3,$
	\item $A_{3}(\mathbf{c})=\left(
	\begin{array}{cccc}
	0 & 1 & 1 & 0 \\
	\beta _1& \beta _2 & 1&-1
	\end{array}\right),\ \mbox{where}\ \mathbf{c}=(\beta_1, \beta_2)\in \mathbb{F}^2,$
	\item $A_{4}(\mathbf{c})=\left(
	\begin{array}{cccc}
	\alpha _1 & 0 & 0 & 0 \\
	0 & \beta _2& 1-\alpha _1&0
	\end{array}\right),\ \mbox{where}\ \mathbf{c}=(\alpha_1, \beta_2)\in \mathbb{F}^2,$
	\item $A_{5}(\mathbf{c})=\left(
	\begin{array}{cccc}
	\alpha_1& 0 & 0 & 0 \\
	1 & 2\alpha_1-1 & 1-\alpha_1&0
	\end{array}\right),\ \mbox{where}\ \mathbf{c}=\alpha_1\in \mathbb{F},$
	\item $A_{6}(\mathbf{c})=\left(
	\begin{array}{cccc}
	\alpha_1 & 0 & 0 & 1 \\
	\beta _1& 1-\alpha_1 & -\alpha_1&0
	\end{array}\right)\simeq \left(
	\begin{array}{cccc}
	\alpha_1 & 0 & 0 & 1 \\
	-\beta _1& 1-\alpha_1 & -\alpha_1&0
	\end{array}\right),\ \mbox{where}\ \mathbf{c}=(\alpha_1, \beta_1)\in \mathbb{F}^2,$
	\item $A_{7}(\mathbf{c})=\left(
	\begin{array}{cccc}
	0 & 1 & 1 & 0 \\
	\beta_1& 1& 0&-1
	\end{array}\right),\ \mbox{where}\ \mathbf{c}=\beta_1\in \mathbb{F},$
	\item $A_{8}(\mathbf{c})=\left(
	\begin{array}{cccc}
	\alpha_1 & 0 & 0 & 0 \\
	0 & 1-\alpha_1 & -\alpha_1&0
	\end{array}\right),\ \mbox{where}\ \mathbf{c}=\alpha_1\in \mathbb{F},$
	\item $A_{9}=\left(
	\begin{array}{cccc}
	\frac{1}{3}& 0 & 0 & 0 \\
	1 & \frac{2}{3} & -\frac{1}{3}&0
	\end{array}\right),\ \ A_{10}=\left(
	\begin{array}{cccc}
	0 & 1 & 1 & 0 \\
	0 &0& 0 &-1
	\end{array}
	\right),\ \ A_{11}=\left(
	\begin{array}{cccc}
	0 & 1 & 1 & 0 \\
	1 &0& 0 &-1
	\end{array}
	\right),\\ A_{12}=\left(
	\begin{array}{cccc}
	0 & 0 & 0 & 0 \\
	1 &0&0 &0\end{array}
	\right).$
\end{itemize}\end{theorem}

\begin{theorem}\label{thm2} Over an algebraically closed field $\mathbb{F}$ $(Char(\mathbb{F})=2)$, any non-trivial $2$-dimensional algebra is isomorphic to only one of the following algebras listed by their matrices of structure constants:
\begin{itemize}
\item $A_{1,2}(\mathbf{c})=\left(
\begin{array}{cccc}
\alpha_1 & \alpha_2 &\alpha_2+1 & \alpha_4 \\
\beta_1 & -\alpha_1 & -\alpha_1+1 & -\alpha_2
\end{array}\right),\ \mbox{where}\ \mathbf{c}=(\alpha_1, \alpha_2, \alpha_4, \beta_1)\in \mathbb{F}^4,$
\item $A_{2,2}(\mathbf{c})=\left(
\begin{array}{cccc}
\alpha_1 & 0 & 0 & 1 \\
\beta _1& \beta_2 & 1-\alpha_1&0
\end{array}\right),\ \mbox{where}\ \mathbf{c}=(\alpha_1, \beta_1, \beta_2)\in \mathbb{F}^3,$
\item $A_{3,2}(\mathbf{c})=\left(
\begin{array}{cccc}
\alpha_1 & 1 & 1 & 0 \\
0& \beta_2 & 1-\alpha_1&1
\end{array}\right),\ \mbox{where}\ \mathbf{c}=(\alpha_1, \beta_2)\in \mathbb{F}^2,$
\item $A_{4,2}(\mathbf{c})=\left(
\begin{array}{cccc}
\alpha _1 & 0 & 0 & 0 \\
0 & \beta_2 & 1-\alpha _1&0
\end{array}\right),\ \mbox{where}\ \mathbf{c}=(\alpha_1,\beta_2)\in \mathbb{F}^2,$
\item $A_{5,2}(\mathbf{c})=\left(
\begin{array}{cccc}
\alpha_1 & 0 & 0 & 0 \\
1 & 1 & 1-\alpha_1&0
\end{array}\right),\ \mbox{where}\ \mathbf{c}=\alpha_1\in \mathbb{F},$
\item $A_{6,2}(\mathbf{c})=\left(
\begin{array}{cccc}
\alpha_1 & 0 & 0 & 1 \\
\beta _1& 1-\alpha_1 & -\alpha_1&0
\end{array}\right),\ \mbox{where}\ \mathbf{c}=(\alpha_1, \beta_1)\in \mathbb{F}^2,$
\item $A_{7,2}(\mathbf{c})=\left(
\begin{array}{cccc}
\alpha_1 & 1 & 1 & 0 \\
0& 1-\alpha_1& -\alpha_1&-1
\end{array}\right),\ \mbox{where}\ \mathbf{c}=\alpha_1\in \mathbb{F},$
\item $A_{8,2}(\mathbf{c})=\left(
\begin{array}{cccc}
\alpha_1 & 0 & 0 & 0 \\
0 & 1-\alpha_1 & -\alpha_1&0
\end{array}\right),\ \mbox{where}\ \mathbf{c}=\alpha_1\in \mathbb{F},$
\item $A_{9,2}=\left(
\begin{array}{cccc}
1 & 0 & 0 & 0 \\
1 & 0 & 1&0
\end{array}\right),\ \ A_{10,2}=\left(
\begin{array}{cccc}
0 & 1 & 1 & 0 \\
0 &0& 0 &-1
\end{array}
\right),\ \ A_{11,2}=\left(
\begin{array}{cccc}
1 & 1 & 1 & 0 \\
0 &-1& -1 &-1
\end{array}
\right),\\ A_{12,2}=\left(
\begin{array}{cccc}
0 & 0 & 0 & 0 \\
1 &0&0 &0\end{array}
\right).$
\end{itemize}
\end{theorem}

\begin{theorem}\label{thm3} Over an algebraically closed field $\mathbb{F}$ $(Char(\mathbb{F})=3)$, any non-trivial $2$-dimensional algebra is isomorphic to only one of the following algebras listed by their matrices of structure constant matrices:
\begin{itemize}
\item $A_{1,3}(\mathbf{c})=\left(
\begin{array}{cccc}
\alpha_1 & \alpha_2 &\alpha_2+1 & \alpha_4 \\
\beta_1 & -\alpha_1 & -\alpha_1+1 & -\alpha_2
\end{array}\right),\ \mbox{where}\ \mathbf{c}=(\alpha_1, \alpha_2, \alpha_4, \beta_1)\in \mathbb{F}^4,$
\item $A_{2,3}(\mathbf{c})=\left(
\begin{array}{cccc}
\alpha_1 & 0 & 0 & 1 \\
\beta _1& \beta _2& 1-\alpha_1&0
\end{array}\right)\simeq\left(
\begin{array}{cccc}
\alpha_1 & 0 & 0 & 1 \\
-\beta _1& \beta _2& 1-\alpha_1&0
\end{array}\right),\ \mbox{where}\ \mathbf{c}=(\alpha_1, \beta_1, \beta_2)\in \mathbb{F}^3,$
\item $A_{3,3}(\mathbf{c})=\left(
\begin{array}{cccc}
0 & 1 & 1 & 0 \\
\beta _1& \beta _2 & 1&-1
\end{array}\right),\ \mbox{where}\ \mathbf{c}=(\beta_1, \beta_2)\in \mathbb{F}^2,$
\item $A_{4,3}(\mathbf{c})=\left(
\begin{array}{cccc}
\alpha _1 & 0 & 0 & 0 \\
0 & \beta _2& 1-\alpha _1&0
\end{array}\right),\ \mbox{where}\ \mathbf{c}=(\alpha_1, \beta_2)\in \mathbb{F}^2,$
\item $A_{5,3}(\mathbf{c})=\left(
\begin{array}{cccc}
\alpha_1& 0 & 0 & 0 \\
1 & -1-\alpha_1 & 1-\alpha_1&0
\end{array}\right),\ \mbox{where}\ \mathbf{c}=\alpha_1\in \mathbb{F},$
\item $A_{6,3}(\mathbf{c})=\left(
\begin{array}{cccc}
\alpha_1 & 0 & 0 & 1 \\
\beta _1& 1-\alpha_1 & -\alpha_1&0
\end{array}\right)\simeq \left(
\begin{array}{cccc}
\alpha_1 & 0 & 0 & 1 \\
-\beta _1& 1-\alpha_1 & -\alpha_1&0
\end{array}\right),\ \mbox{where}\ \mathbf{c}=(\alpha_1, \beta_1)\in \mathbb{F}^2,$
\item $A_{7,3}(\mathbf{c})=\left(
\begin{array}{cccc}
0 & 1 & 1 & 0 \\
\beta_1& 1& 0&-1
\end{array}\right),\ \mbox{where}\ \mathbf{c}=\beta_1\in \mathbb{F},$
\item $A_{8,3}(\mathbf{c})=\left(
\begin{array}{cccc}
\alpha_1 & 0 & 0 & 0 \\
0 & 1-\alpha_1 & -\alpha_1&0
\end{array}\right),\ \mbox{where}\ \mathbf{c}=\alpha_1\in \mathbb{F},$
\item $A_{9,3}=\left(
\begin{array}{cccc}
0 & 1& 1& 0 \\
1 &0&0 &-1\end{array}
\right),\ \ A_{10,3}=\left(
\begin{array}{cccc}
0 & 1 & 1 & 0 \\
0 &0&0 &-1\end{array}
\right),\ \ A_{11,3}=\left(
\begin{array}{cccc}
1 & 0 & 0 & 0 \\
1 &-1&-1 &0\end{array}
\right),\\ A_{12,3}=\left(
\begin{array}{cccc}
0 &0 &0 & 0 \\
1 &0& 0 &0
\end{array}
\right).$
\end{itemize}\end{theorem}
\begin{remark} In reality in \cite{A} the class $A_{3,2}(\mathbf{c})$ should be understood as it is in this paper, as far as there is a type-mistake in this case in \cite{A}.\end{remark}

 \section{Two-dimensional left unital algebras}

Let $\mathbb{A}$ be a left unital algebra. In terms of its MSC $A$ the algebra $\mathbb{A}$ to be left unital is written as follows:
\begin{equation}\label{7}
  A(l\otimes u)=u,
\end{equation} where $u = (u_1, u_2, ..., u_m)^T,$ and $l = (t_1, t_2, ..., t_m)^T$ are column coordinate vectors of $\mathbf{u}$ and $\mathbf{1_\textit{L}},$ respectively.\\

It is easy to see that for a given $2$-dimensional algebra $\mathbb{A}$ with MSC $A=\left(\begin{array}{cccc} \alpha_1 & \alpha_2 & \alpha_3 &\alpha_4\\ \beta_1
& \beta_2 & \beta_3 &\beta_4\end{array}\right)$ the existence of a left unit element is equivalent to the equality of ranks of the matrices  \[M=\left(\begin{array}{cc} \alpha_1 & \alpha_3\\
\beta_1 & \beta_3\\ \alpha_2 & \alpha_4\\ \beta_2-\alpha_1 & \beta_4-\alpha_3
\end{array}\right)\ \mbox{and}\ M'= \left(\begin{array}{ccc} \alpha_1 & \alpha_3&1\\
\beta_1 & \beta_3&0\\ \alpha_2 & \alpha_4&0\\ \beta_2-\alpha_1 & \beta_4-\alpha_3&0
\end{array}\right).\] This equality holds if and only if
\begin{equation}\label{1}  \left|\begin{array}{cc}\beta_1&\beta_3\\ \alpha_2&\alpha_4\end{array}\right|=\left|\begin{array}{cc}\beta_1&\beta_3\\ \beta_2-\alpha_1&\beta_4-\alpha_3\end{array}\right|=\left|\begin{array}{cc}\alpha_2&\alpha_4\\ \beta_2-\alpha_1&\beta_4-\alpha_3\end{array}\right|=0,\end{equation}  and at least one of the following two cases holds true:
\begin{equation}\label{2} (\alpha_1,\alpha_3)\neq 0,\  (\beta_1,\beta_3)=(\alpha_2,\alpha_4)=(\beta_2-\alpha_1,\beta_4-\alpha_3)=0,\end{equation} or \begin{equation}\label{3} \left|\begin{array}{cc}\alpha_1&\alpha_3\\ a&b\end{array}\right|\neq 0,\ \ \mbox{whenever there exists nonzero}\ (a,b)\in \{(\beta_1,\beta_3),(\alpha_2,\alpha_4),(\beta_2-\alpha_1,\beta_4-\alpha_3)\}. \end{equation} 

Note that the conditions (\ref{1}), (\ref{2}) and (\ref{1}), (\ref{3}) correspond to the existence of many and unique left units, respectively.
\begin{theorem}\label{thm4} Over any algebraically closed field $\mathbb{F}$ $\left(Char\left(\mathbb{F}\right)\neq 2\ \mbox{and}\ 3\right)$ any nontrivial $2$-dimensional left unital algebra is isomorphic to only one of the following non-isomorphic left unital algebras presented by their MSC:
\begin{itemize}
\item $A_{1}\left(\alpha_1,\frac{\alpha_1(1-\alpha_1)}{\beta_1}-\frac{1}{2},\frac{\alpha_1(1-\alpha_1)^2}{\beta^2_1}-\frac{1-\alpha_1}{2\beta_1},\beta_1\right)
\\ =\left(
\begin{array}{cccc}
 \alpha _1 & \frac{2 \alpha _1-2 \alpha _1^2-\beta _1}{2 \beta _1} & \frac{2 \alpha _1-2 \alpha _1^2+\beta _1}{2 \beta _1} & \frac{2 \alpha _1-4 \alpha _1^2+2 \alpha _1^3-\beta _1+\alpha _1 \beta _1}{2 \beta _1^2} \\
 \beta _1 & -\alpha _1 & 1-\alpha _1 & \frac{-2 \alpha _1+2 \alpha _1^2+\beta _1}{2 \beta _1}
\end{array}
\right),$ where $\beta_1\neq 0,$
\item $A_{1}\left(1,\alpha_2,\frac{\alpha_2(2\alpha_2+1)}{2},0\right)=\left(
\begin{array}{cccc}
 1 & \alpha _2 & 1+\alpha _2 & \frac{1}{2} \left(\alpha _2+2 \alpha _2^2\right) \\
 0 & -1 & 0 & -\alpha _2
\end{array}
\right),$
\item $A_{2}(\alpha_1,0,\alpha_1)=\left(
\begin{array}{cccc}
 \alpha _1 & 0 & 0 & 1 \\
 0 & \alpha _1 & -\alpha _1+1 & 0
\end{array}
\right),$ where $\alpha _1\neq 0,$
\item $A_{4}(\alpha_1,\alpha_1)=\left(
\begin{array}{cccc}
 \alpha _1 & 0 & 0 & 0 \\
 0 & \alpha _1 & -\alpha _1+1 & 0
\end{array}
\right),$ where $\alpha _1\neq 0,$
\item $A_{6}\left(\frac{1}{2},0\right)=\left(
\begin{array}{cccc}
 \frac{1}{2} & 0 & 0 & 1 \\
 0 & \frac{1}{2} & -\frac{1}{2} & 0
\end{array}
\right),\ \ A_{8}\left(\frac{1}{2}\right)=\left(
\begin{array}{cccc}
 \frac{1}{2} & 0 & 0 & 0 \\
 0 & \frac{1}{2} & -\frac{1}{2} & 0
\end{array}
\right).$
  \end{itemize}
\end{theorem}
\textbf{Proof.} Let us consider $A_{1}(\mathbf{c})=\left(
\begin{array}{cccc}
\alpha_1 & \alpha_2 &\alpha_2+1 & \alpha_4 \\
\beta_1 & -\alpha_1 & -\alpha_1+1 & -\alpha_2
\end{array}\right)$.\\ Then $M=\left(
\begin{array}{cc}
\alpha_1 & \alpha_2+1\\ \beta_1 &1-\alpha_1\\ \alpha_2& \alpha_4 \\
-2\alpha_1 & -2\alpha_2-1\end{array}\right)$ and the equality
 (\ref{1}) means
\[\beta_1\alpha_4-\alpha_2(1-\alpha_1)=-\beta_1(2\alpha_2+1)+2\alpha_1(1-\alpha_1)=-\alpha_2(2\alpha_2+1)+2\alpha_1\alpha_4=0\]
and  (\ref{2}) doesn't occur. There are two possibilities:\\
\underline{\textbf{Case 1.} $\beta_1\neq 0.$} In this case the equality (\ref{1}) is equivalent to\\
$\alpha_4=\frac{\alpha_2(1-\alpha_1)}{\beta_1}, \alpha_2=\frac{\alpha_1(1-\alpha_1)}{\beta_1}-\frac{1}{2}, \
\mbox{and}\ 
\left|\begin{array}{cc}
\alpha_1 & \alpha_2+1\\\beta_1 & 1-\alpha_1 \end{array}\right|=-\beta_1/2\neq 0$.\\ 
Therefore, $A_1\left(\alpha_1,\frac{\alpha_1(1-\alpha_1)}{\beta_1}-\frac{1}{2},\frac{\alpha_1(1-\alpha_1)^2}{\beta^2_1}-\frac{1-\alpha_1}{2\beta_1},\beta_1\right)$  has a left unit, where $\beta_1\neq 0$. \\
\underline{\textbf{Case 2.} $\beta_1=0.$} In this case the equality (\ref{1}) is equivalent to \[\alpha_2(1-\alpha_1)=\alpha_1(1-\alpha_1)=-\alpha_2(2\alpha_2+1)+2\alpha_1\alpha_4=0\] and 
 (\ref{3}) occurs if and only if $\alpha_1=1$ and therefore \[A_1\left(1,\alpha_2,\frac{\alpha_2(2\alpha_2+1)}{2},0\right)\] also has a left unit.

Consider $A_{2}(\mathbf{c})=\left(
\begin{array}{cccc}
\alpha_1 & 0 & 0 & 1 \\
\beta _1& \beta _2& 1-\alpha_1&0
\end{array}\right)$. Then $M=\left(
\begin{array}{cc}
\alpha_1 & 0\\ \beta_1 &1-\alpha_1\\ 0& 1 \\
\beta_2-\alpha_1 & 0\end{array}\right).$\\ The equality
 (\ref{1}) means \[\beta_1=(1-\alpha_1)(\beta_2-\alpha_1)=\beta_2-\alpha_1= 0\] and  (\ref{2}) doesn't occur. Therefore, $A_2(\alpha_1,0,\alpha_1)$ has a left unit, where $\alpha _1\neq 0$.

In $A_{3}(\mathbf{c})=\left(
\begin{array}{cccc}
0 & 1 & 1 & 0 \\
\beta _1& \beta _2 & 1&-1
\end{array}\right)$ case we have $M=\left(
\begin{array}{cc}
\alpha_1 & 1\\ \beta_1 &1\\ 1& 0 \\
\beta_2 & -2\end{array}\right)$ and $\left|\begin{array}{cc}1& 0\\ \beta_2&-2\end{array}\right|=-2\neq0,$ which shows the absence of a left unit.

Let us consider $A_{4}(\mathbf{c})=\left(
\begin{array}{cccc}
\alpha _1 & 0 & 0 & 0 \\
0 & \beta _2& 1-\alpha _1&0
\end{array}\right)$. Then $M=\left(
\begin{array}{cc}
\alpha_1 & 0\\ 0 &1-\alpha_1\\ 0& 0 \\
\beta_2-\alpha_1 & 0\end{array}\right)$, the equality (\ref{1}) is equivalent to $(1-\alpha_1)(\alpha _1-\beta_2)=0$ and therefore $A_4(1, 1)$ has  left units. In this case  (\ref{3}) happens if and only if  $\alpha _1\neq 0,1$, $\alpha _1=\beta_2$. So $A_{4}(\alpha_1, \alpha_1)$ has a left unit, where $\alpha_1\neq 0.$

In $A_{5}(\mathbf{c})=\left(
\begin{array}{cccc}
\alpha_1& 0 & 0 & 0 \\
1 & 2\alpha_1-1 & 1-\alpha_1&0
\end{array}\right)$ case one has $M=\left(
\begin{array}{cc}
\alpha_1 & 0\\ 1 &1-\alpha_1\\ 0& 0 \\
\alpha_1-1 & 0\end{array}\right)$, the equality (\ref{1}) means $(1-\alpha _1)(\alpha_1-1)=0$, so $\alpha_1=1.$ But neither (\ref{2}) no (\ref{3}) occurs, that is among $A_5(\alpha_1)$ there is no algebra with a left unit.

In $A_{6}(\mathbf{c})=\left(
\begin{array}{cccc}
\alpha_1 & 0 & 0 & 1 \\
\beta _1& 1-\alpha_1 & -\alpha_1&0
\end{array}\right)$ case we have $M=\left(
\begin{array}{cc}
\alpha_1 & 0\\ \beta_1 &-\alpha_1\\ 0& 1 \\
1-2\alpha_1 & 0\end{array}\right)$, the equality (\ref{1}) is equivalent to $\beta_1=\alpha_1(1-2\alpha_1)= -1+2\alpha_1=0$ and therefore $A_{6}\left(\frac{1}{2},0\right)$ has a left unit.

In $A_{7}(\mathbf{c})=\left(
\begin{array}{cccc}
0 & 1 & 1 & 0 \\
\beta_1& 1& 0&-1
\end{array}\right)$ case we have $M=\left(
\begin{array}{cc}
0 & 1\\ \beta_1 &0\\ 1& 0 \\
1 & -2\end{array}\right)$, and the inequality\\ $\left|\begin{array}{cc}1& 0\\ 1&-2\end{array}\right|=-2\neq0$ shows the absence of a left unit due to (\ref{1}).

In $A_{8}(\mathbf{c})=\left(
\begin{array}{cccc}
\alpha_1 & 0 & 0 & 0 \\
0 & 1-\alpha_1 & -\alpha_1&0
\end{array}\right)$ case $M=\left(
\begin{array}{cc}
\alpha_1 & 0\\ 0 &-\alpha_1\\ 0& 0 \\
1-2\alpha_1 & 0\end{array}\right)$, the equality (\ref{1}) gives\\ $\alpha_1(1-2\alpha_1)=0$ and therefore $A_{8}(\frac{1}{2})$ has a left unit.

It is easy to see that for $A_{9}, A_{10}$, $A_{11}$ the equality (\ref{1}) does not occur, the equalities (\ref{2}), (\ref{3}) don't occur for $A_{12}$ and therefore they have no left units.

We present the corresponding results in characteristic $2$ and $3$ cases without proof as follows.

\begin{theorem}\label{5} Over any algebraically closed field $\mathbb{F}$ of characteristic $2$ any nontrivial $2$-dimensional left unital algebra is isomorphic to only one of the following non-isomorphic left unital algebras  presented by their MSC:
	\begin{itemize}
		\item $A_{1,2}\left(\alpha_1,0,\alpha_4,0\right)=\left(
		\begin{array}{cccc}
		\alpha _1 & 0 & 1 & \alpha _4 \\
		0 & -\alpha _1 & 1-\alpha _1 & 0
		\end{array}
		\right),$ where $\alpha_1\neq 0,$
		\item $A_{2,2}(\alpha_1,0,\alpha_1)=\left(
		\begin{array}{cccc}
		\alpha _1 & 0 & 0 & 1 \\
		0 & \alpha _1 & -\alpha _1+1 & 0
		\end{array}
		\right),$ where $\alpha _1\neq 0,$
		\item $A_{3,2}(1,\beta_2)=\left(
		\begin{array}{cccc}
		1 & 1 & 1 & 0 \\
		0 & \beta_2 & 0 & 1	\end{array}
		\right),$ 
				\item $A_{4,2}(\alpha_1,\alpha_1)=\left(
		\begin{array}{cccc}
		\alpha _1 & 0 & 0 & 0 \\
		0 & \alpha _1 & -\alpha _1+1 & 0
		\end{array}
		\right),$ where $\alpha _1\neq 0,$
		\item $A_{7,2}\left(0\right)=\left(
		\begin{array}{cccc}
		0 & 1 & 1 & 0 \\
		0 & 1 & 0 & -1
		\end{array}
		\right),\ \ A_{10,2}=\left(
		\begin{array}{cccc}
		0 & 1 & 1 & 0 \\
		0 & 0 & 0 & 1
		\end{array}
		\right).$
	\end{itemize}
\end{theorem}
\begin{theorem} Over any algebraically closed field $\mathbb{F}$ of characteristic $3$ any nontrivial $2$-dimensional left unital algebra is isomorphic to only one of the following non-isomorphic left unital algebras  presented by their MSC:
	\begin{itemize}
		\item $A_{1,3}\left(\alpha_1,\frac{\alpha_1(1-\alpha_1)}{\beta_1}-\frac{1}{2},\frac{\alpha_1(1-\alpha_1)^2}{\beta^2_1}-\frac{1-\alpha_1}{2\beta_1},\beta_1\right)
		\\ =\left(
		\begin{array}{cccc}
		\alpha _1 & \frac{2 \alpha _1-2 \alpha _1^2-\beta _1}{2 \beta _1} & \frac{2 \alpha _1-2 \alpha _1^2+\beta _1}{2 \beta _1} & \frac{2 \alpha _1-4 \alpha _1^2+2 \alpha _1^3-\beta _1+\alpha _1 \beta _1}{2 \beta _1^2} \\
		\beta _1 & -\alpha _1 & 1-\alpha _1 & \frac{-2 \alpha _1+2 \alpha _1^2+\beta _1}{2 \beta _1}
		\end{array}
		\right),$ where $\beta_1\neq 0,$
		\item $A_{1,3}\left(1,\alpha_2,\frac{\alpha_2(2\alpha_2+1)}{2},0\right)=\left(
		\begin{array}{cccc}
		1 & \alpha _2 & 1+\alpha _2 & \frac{1}{2} \left(\alpha _2+2 \alpha _2^2\right) \\
		0 & -1 & 0 & -\alpha _2
		\end{array}
		\right),$
		\item $A_{2,3}(\alpha_1,0,\alpha_1)=\left(
		\begin{array}{cccc}
		\alpha _1 & 0 & 0 & 1 \\
		0 & \alpha _1 & -\alpha _1+1 & 0
		\end{array}
		\right),$ where $\alpha _1\neq 0,$
		\item $A_{4,3}(\alpha_1,\alpha_1)=\left(
		\begin{array}{cccc}
		\alpha _1 & 0 & 0 & 0 \\
		0 & \alpha _1 & -\alpha _1+1 & 0
		\end{array}
		\right),$ where $\alpha _1\neq 0,$
		\item $A_{6,3}\left(\frac{1}{2},0\right)=\left(
		\begin{array}{cccc}
		\frac{1}{2} & 0 & 0 & 1 \\
		0 & \frac{1}{2} & -\frac{1}{2} & 0
		\end{array}
		\right),\ \ A_{8,3}\left(\frac{1}{2}\right)=\left(
		\begin{array}{cccc}
		\frac{1}{2} & 0 & 0 & 0 \\
		0 & \frac{1}{2} & -\frac{1}{2} & 0
		\end{array}
		\right).$
	\end{itemize}
\end{theorem}
Note that according to Theorem 3.1 and Theorem 3.3 in the cases $Char(\mathbb{F})\neq 2$,$3$ and $Char(\mathbb{F})=3$  the lists are identical. Therefore, we summarize the final result for $2$-dimensional left unital algebras in the following table, where all left units as well are given. 
\begin{table}[h]
{	\begin{tabular}{ |c|c|c| }
		\hline
		&$\mbox{Algebra}$ & $\mathbf{1_\textit{L}} $ \\
		\hline
		\multirow{7}{*}{\begin{turn}{90}$Char\left(\mathbb{F}\right)\neq 2 \ \ \ \ \ \ $\end{turn}}
		&$A_{1}\left(\alpha_1,\frac{\alpha_1(1-\alpha_1)}{\beta_1}-\frac{1}{2},\frac{\alpha_1(1-\alpha_1)^2}{\beta^2_1}-\frac{1-\alpha_1}{2\beta_1},\beta_1\right)$, where $\beta_1\neq 0$&  $
		\frac{-2(1-\alpha_1)}{\beta_1}e_1+2e_2 
				$
		\\
		\cline{2-3}
		&$A_1\left(1,\alpha_2,\frac{\alpha_2(2\alpha_2+1)}{2},0\right)$ &  $
		-(1+2\alpha_2)e_1+2e_2
		$
		\\
		\cline{2-3}
		&$A_2(\alpha_1,0,\alpha_1)$, where $\alpha_1\neq 0$&$
		\frac{1}{\alpha_1}e_1
		$
		\\
		\cline{2-3}
		&$A_4(1, 1)$&$e_1+te_2
		,$ where $t \in \mathbb{F}$
		\\
		\cline{2-3}
		&$A_{4}(\alpha_1, \alpha_1)$, where $\alpha_1\neq 0,1$
		&$
		\frac{1}{\alpha_1}e_1
		$
		\\
		\cline{2-3}
		&$A_{6}\left(\frac{1}{2},0\right)$&$2e_1
		$
		\\
		\cline{2-3}
		&$A_{8}(\frac{1}{2})$&$2e_1
		$
		\\
		\hline
		\multirow{6}{*}{\begin{turn}{90}$Char\left(\mathbb{F}\right)= 2 \ \ \ \ \ \  $\end{turn}}&$A_{1,2}\left(\alpha_1,0,\alpha_4,0\right),$ where $\alpha_1\neq 0$& $
		\frac{1}{\alpha_1}e_1
		$
		\\
		\cline{2-3}
		&$A_{2,2}(\alpha_1,0,\alpha_1),$ where $\alpha _1\neq 0$ &$
		\frac{1}{\alpha_1}e_1
		$
		\\
		\cline{2-3}
		&$A_{3,2}(1,\beta_2)$  &$e_2
		$
		\\
		\cline{2-3}
		
		&$A_{4,2}(1,1)$&$
		e_1+te_2
		$, where $t \in \mathbb{F}$
		\\
		\cline{2-3}
		&$A_{4,2}(\alpha_1,\alpha_1),$ where $\alpha _1\neq 0,1$ &$
		\frac{1}{\alpha_1}e_1
		$
		\\
		\cline{2-3}
		&$A_{7,2}\left(0\right)$&$e_2
		$
		\\
		\cline{2-3}
		&$A_{10,2}$&$e_2
		$
		\\
		\hline
	\end{tabular}}\end{table}

  \section{Two-dimensional right unital algebras}
Now let us consider the existence of a right unit for an algebra $\mathbb{A}$ given by its MSC $A=$\\  $\left(\begin{array}{cccc} \alpha_1 & \alpha_2 & \alpha_3 &\alpha_4\\ \beta_1 & \beta_2 & \beta_3 &\beta_4\end{array}\right)$. It is easy to see that $\mathbb{A}$ has a right unit element if and only if the following matrices \[\left(\begin{array}{cc} \alpha_1 & \alpha_2\\
\beta_1 & \beta_2\\ \alpha_3 & \alpha_4\\ \beta_3-\alpha_1 & \beta_4-\alpha_2
\end{array}\right), \left(\begin{array}{ccc} \alpha_1 & \alpha_2&1\\
\beta_1 & \beta_2&0\\ \alpha_3 & \alpha_4&0\\ \beta_3-\alpha_1 & \beta_4-\alpha_2&0
\end{array}\right)\] have equal ranks. It happens if and only if
\[\left|\begin{array}{cc}\beta_1&\beta_2\\ \alpha_3&\alpha_4\end{array}\right|=\left|\begin{array}{cc}\beta_1&\beta_2\\ \beta_3-\alpha_1&\beta_4-\alpha_2\end{array}\right|=\left|\begin{array}{cc}\alpha_3&\alpha_4\\ \beta_3-\alpha_1&\beta_4-\alpha_2\end{array}\right|=0
\] and at least one of the following two cases holds true
\[(\alpha_1,\alpha_2)\neq 0,\  (\beta_1,\beta_2)=(\alpha_3,\alpha_4)=(\beta_3-\alpha_1,\beta_4-\alpha_2)=0,
                    \] or
\[\left|\begin{array}{cc}\alpha_1&\alpha_2\\ a&b\end{array}\right|\neq 0,\ \ \mbox{if there exists nonzero}\ (a,b)\in \{(\beta_1,\beta_2),(\alpha_3,\alpha_4),(\beta_3-\alpha_1,\beta_4-\alpha_2)\}.\]

Because of similarity of proofs in right unital cases to those of left unital ones we present the result without proof by the following theorems.

\begin{theorem}\label{thm7} Over any algebraically closed field $\mathbb{F}$ of characteristic not $2$ any nontrivial $2$-dimensional right unital algebra is isomorphic to only one of the following non-isomorphic right unital algebras presented by their MSC:
\begin{itemize}
\item $A_{1}\left(\alpha_1,\frac{\alpha_1(1-2\alpha_1)}{2\beta_1},-\frac{\alpha^2_1(1-2\alpha_1)}{2\beta^2_1}-\frac{\alpha_1}{\beta_1},\beta_1\right),$ where $\alpha_1\beta_1\neq 0,$
	\item $A_{1}(0,\alpha_2,-2\alpha_2(1+\alpha_2),0),$ where $\alpha_2(1+\alpha_2)\neq 0,$
 \item $A_{1}\left(\frac{1}{2},-1,\alpha_4,0\right),\ \ A_{2}\left(\frac{1}{2},0,\beta_2\right),\ \ A_{4}\left(\frac{1}{2},\beta_2\right).$
\end{itemize}
\end{theorem}

\begin{theorem}\label{thm6} Over any algebraically closed field $\mathbb{F}$ of characteristic $2$ any nontrivial $2$-dimensional right unital algebra is isomorphic to only one of the following non-isomorphic right unital algebras presented by their MSC:
\begin{itemize}
  \item $A_{1,2}\left(0,\alpha_2,0,\beta_1\right),$ where $\alpha_2\neq 0,$
  \item $A_{3,2}(\alpha_1,0),$
  \item $A_{6,2}\left(\alpha_1,0\right),$ where $\alpha_1\neq 0,$
  \item $A_{7,2}\left(1\right),$
  \item $A_{8,2}(\alpha_1),$ where $\alpha_1\neq 0,$
  \item $A_{10,2}.$
\end{itemize}
\end{theorem}

The results obtained are summarized in the following table, where all right units as well are listed. 
\begin{table}[h]
{
	\begin{tabular}{ |c|c|c| }
		\hline
		&$\mbox{Algebra}$ & $\mathbf{1_\textit{R}} $ \\
		\hline
		\multirow{6}{*}{\begin{turn}{90}$Char\left(\mathbb{F}\right)\neq 2 \ \ \ \ \  $\end{turn}}&$A_{1}\left(\alpha_1,\frac{\alpha_1(1-2\alpha_1)}{2\beta_1},\frac{-\alpha^2_1(1-2\alpha_1)}{2\beta^2_1}-\frac{\alpha_1}{\beta_1},\beta_1\right)$, where $\alpha_1\beta_1\neq 0$& $2e_1+
		\frac{2\beta_1}{\alpha_1}e_2 
		$
		\\
		\cline{2-3}
		&$A_{1}(0,\alpha_2, -2\alpha_2(1+\alpha_2),0)$, where $\alpha_2(1+\alpha_2)\neq 0$ & $2e_1+	\frac{1}{\alpha_2}e_2
		$
		\\
		\cline{2-3}
		&$A_{1}\left(\frac{1}{2},-1,\alpha_4,0\right)$ & $
		2e_1
		$
		\\
		\cline{2-3}
		&$A_{2}(\frac{1}{2}, 0,\beta_2)$ &$2e_1
		$
		\\
		\cline{2-3}
		&$A_{4}(\frac{1}{2},0)$&$2e_1+te_2
		,$ where $t \in \mathbb{F}$
		\\
		\cline{2-3}
		&$A_{4}(\frac{1}{2},\beta_2)$, where $\beta_2\neq0$
		&$2e_1
		$
		\\
		\hline
		\multirow{7}{*}{\begin{turn}{90}$Char\left(\mathbb{F}\right)= 2 \ \ \ \ \  $\end{turn}}&$A_{1,2}\left(0,\alpha_2,0,\beta_1\right),$ where $\alpha_2\neq 0$& $
		\frac{1}{\alpha_2}e_2
		$
		\\
		\cline{2-3}
		&$A_{3,2}(\alpha_1,0)$&$e_2
		$
		\\
		\cline{2-3}
		&$A_{6,2}(\alpha_1,0),$ where $\alpha _1\neq 0$&$
		\frac{1}{\alpha_1}e_1
		$
		\\
		\cline{2-3}
		&$A_{7,2}\left(1\right)$&$e_2
		$
		\\
		\cline{2-3}
		&$A_{8,2}\left(1\right)$&$e_1+te_2
		$, where $t \in \mathbb{F}$
		\\
		\cline{2-3}
		&$A_{8,2}\left(\alpha_1\right)$, where $\alpha_1 \neq 0, 1$ &$e_1
		$
		\\
		\cline{2-3}
		&$A_{10,2}$&$e_2
		$
		\\
		\hline
\end{tabular}}\end{table}

\begin{corollary} 
	Over an algebraically closed field $\mathbb{F}$, $\left(Char(\mathbb{F})\neq 2\right)$, there exist, up to isomorphism, only two non-trivial $2$-dimensional
unital algebras given by their matrices of structure constants as follows
\[A_{2}\left(\frac{1}{2},0,\frac{1}{2}\right)=\left(
	\begin{array}{cccc}
	\frac{1}{2} & 0 &0 & 1 \\
	0 & \frac{1}{2} & \frac{1}{2} & 0
	\end{array}\right),\ \ A_{4}\left(\frac{1}{2},\frac{1}{2}\right)=\left(
	\begin{array}{cccc}
	\frac{1}{2} & 0 &0 & 0 \\
	0 & \frac{1}{2} & \frac{1}{2} & 0
	\end{array}\right).\]\end{corollary}
\begin{corollary} Over an algebraically closed field $\mathbb{F}$, $\left(Char\left(\mathbb{F}\right)=2\right)$, there exists, up to isomorphism, only one non-trivial $2$-dimensional
unital algebra given by its matrix of structure constants as \[ A_{10,2}=\left(
\begin{array}{cccc}
0 & 1 & 1 & 0 \\
0 &0& 0 &1
\end{array}
\right).\]\end{corollary}
\section{Two-dimensional left and right unital real algebras}
Due to \cite{B1} we have the following classification theorem.
\begin{theorem}\label{th1} Any non-trivial 2-dimensional real algebra is isomorphic to only one of the following listed, by their matrices of structure constants, algebras:\\
		$A_{1,r}(\mathbf{c})=\left(
	\begin{array}{cccc}
	\alpha_1 & \alpha_2 &\alpha_2+1 & \alpha_4 \\
	\beta_1 & -\alpha_1 & -\alpha_1+1 & -\alpha_2
	\end{array}\right),\ \mbox{where}\ \mathbf{c}=(\alpha_1, \alpha_2, \alpha_4, \beta_1)\in \mathbb{R}^4,$\\
	$A_{2,r}(\mathbf{c})=\left(
	\begin{array}{cccc}
	\alpha_1 & 0 & 0 & 1 \\
	\beta _1& \beta _2& 1-\alpha_1&0
	\end{array}\right), \mbox{where}\ \beta_1\geq 0,\ \mathbf{c}=(\alpha_1, \beta_1, \beta_2)\in \mathbb{R}^3,$\\
	$A_{3,r}(\mathbf{c})=\left(
	\begin{array}{cccc}
	\alpha_1 & 0 & 0 & -1 \\
	\beta _1& \beta _2& 1-\alpha_1&0
	\end{array}\right), \mbox{where}\ \beta_1\geq 0,\ \mathbf{c}=(\alpha_1, \beta_1, \beta_2)\in \mathbb{R}^3,$\\
	$A_{4,r}(\mathbf{c})=\left(
	\begin{array}{cccc}
	0 & 1 & 1 & 0 \\
	\beta _1& \beta _2 & 1&-1
	\end{array}\right),\ \mbox{where}\ \mathbf{c}=(\beta_1, \beta_2)\in \mathbb{R}^2,$\\
	$A_{5,r}(\mathbf{c})=\left(
	\begin{array}{cccc}
	\alpha _1 & 0 & 0 & 0 \\
	0 & \beta _2& 1-\alpha _1&0
	\end{array}\right),\ \mbox{where}\ \mathbf{c}=(\alpha_1, \beta_2)\in \mathbb{R}^2,$\\
	$A_{6,r}(\mathbf{c})=\left(
	\begin{array}{cccc}
	\alpha_1& 0 & 0 & 0 \\
	1 & 2\alpha_1-1 & 1-\alpha_1&0
	\end{array}\right),\ \mbox{where}\ \mathbf{c}=\alpha_1\in \mathbb{R},$\\
	$A_{7,r}(\mathbf{c})=\left(
	\begin{array}{cccc}
	\alpha_1 & 0 & 0 & 1 \\
	\beta _1& 1-\alpha_1 & -\alpha_1&0
	\end{array}\right), \mbox{where}\ \beta_1\geq 0,\ \mathbf{c}=(\alpha_1, \beta_1)\in \mathbb{R}^2,$\\
	$A_{8,r}(\mathbf{c})=\left(
	\begin{array}{cccc}
	\alpha_1 & 0 & 0 & -1 \\
	\beta _1& 1-\alpha_1 & -\alpha_1&0
	\end{array}\right), \mbox{where}\ \beta_1\geq 0,\ \mathbf{c}=(\alpha_1, \beta_1)\in \mathbb{R}^2,$\\
	$A_{9,r}(\mathbf{c})=\left(
	\begin{array}{cccc}
	0 & 1 & 1 & 0 \\
	\beta_1& 1& 0&-1
	\end{array}\right),\ \mbox{where}\ \mathbf{c}=\beta_1\in \mathbb{R},$\\
	$A_{10,r}(\mathbf{c})=\left(
	\begin{array}{cccc}
	\alpha_1 & 0 & 0 & 0 \\
	0 & 1-\alpha_1 & -\alpha_1&0
	\end{array}\right),\ \mbox{where}\ \mathbf{c}=\alpha_1\in\mathbb{R},$\\
	$ A_{11,r}=\left(
	\begin{array}{cccc}
	\frac{1}{3}& 0 & 0 & 0 \\
	1 & \frac{2}{3} & -\frac{1}{3}&0
	\end{array}\right),\
	\ A_{12,r}=\left(
	\begin{array}{cccc}
	0 & 1 & 1 & 0 \\
	1 &0&0 &-1
	\end{array}
	\right),$\\
	$ A_{13,r}=\left(
	\begin{array}{cccc}
	0 & 1 & 1 & 0\\
	-1 &0&0 &-1
	\end{array}
	\right),\
	\ A_{14,r}=\left(
	\begin{array}{cccc}
	0 & 1 & 1 & 0 \\
	0 & 0 & 0 & -1
	\end{array}
	\right),\
	\ A_{15,r}=\left(
	\begin{array}{cccc}
	0 & 0 & 0 & 0 \\
	1 &0& 0 &0
	\end{array}\right).$\end{theorem}
Owing to Theorem 5.1 the following results can be proved. 
\begin{theorem} Over the real field $\mathbb{R}$ up to isomorphism there exist only the following nontrivial non-isomorphic two dimensional left unital algebras
		\begin{itemize}
			\item $A_{1,r}\left(\alpha_1,\frac{\alpha_1(1-\alpha_1)}{\beta_1}-\frac{1}{2},\frac{\alpha_1(1-\alpha_1)^2}{\beta^2_1}-\frac{1-\alpha_1}{2\beta_1},\beta_1\right)
			\\ =\left(
			\begin{array}{cccc}
			\alpha _1 & \frac{2 \alpha _1-2 \alpha _1^2-\beta _1}{2 \beta _1} & \frac{2 \alpha _1-2 \alpha _1^2+\beta _1}{2 \beta _1} & \frac{2 \alpha _1-4 \alpha _1^2+2 \alpha _1^3-\beta _1+\alpha _1 \beta _1}{2 \beta _1^2} \\
			\beta _1 & -\alpha _1 & 1-\alpha _1 & \frac{-2 \alpha _1+2 \alpha _1^2+\beta _1}{2 \beta _1}
			\end{array}
			\right),$ where $\beta_1\neq 0,$
			\item $A_{1,r}\left(1,\alpha_2,\frac{\alpha_2(2\alpha_2+1)}{2},0\right)=\left(
			\begin{array}{cccc}
			1 & \alpha _2 & 1+\alpha _2 & \frac{1}{2} \left(\alpha _2+2 \alpha _2^2\right) \\
			0 & -1 & 0 & -\alpha _2
			\end{array}
			\right),$
			\item $A_{2,r}(\alpha_1,0,\alpha_1)=\left(
			\begin{array}{cccc}
			\alpha _1 & 0 & 0 & 1 \\
			0 & \alpha _1 & -\alpha _1+1 & 0
			\end{array}
			\right),$ where $\alpha _1\neq 0,$
			\item $A_{3,r}(\alpha_1,0,\alpha_1)=\left(
			\begin{array}{cccc}
			\alpha _1 & 0 & 0 & -1 \\
			0 & \alpha _1 & -\alpha _1+1 & 0
			\end{array}
			\right),$ where $\alpha _1\neq 0,$
			\item $A_{5,r}(\alpha_1,\alpha_1)=\left(
			\begin{array}{cccc}
			\alpha _1 & 0 & 0 & 0 \\
			0 & \alpha _1 & -\alpha _1+1 & 0
			\end{array}
			\right),$ where $\alpha _1\neq 0,$ \item $A_{7,r}\left(\frac{1}{2},0\right)=\left(
			\begin{array}{cccc}
			\frac{1}{2} & 0 & 0 & 1 \\
			0 & \frac{1}{2} & -\frac{1}{2} & 0
			\end{array}
			\right),$
			$ A_{8,r}\left(\frac{1}{2},0\right)=\left(
			\begin{array}{cccc}
			\frac{1}{2} & 0 & 0 & -1 \\
			0 & \frac{1}{2} & -\frac{1}{2} & 0
			\end{array}
			\right),\\  A_{10,r}\left(\frac{1}{2}\right)=\left(
			\begin{array}{cccc}
			\frac{1}{2} & 0 & 0 & 0 \\
			0 & \frac{1}{2} & -\frac{1}{2} & 0
			\end{array}
			\right).$
		\end{itemize}
\end{theorem}
\begin{theorem} Over the real field $\mathbb{R}$ up to isomorphism there exist only the following nontrivial non-isomorphic two dimensional right unital algebras:
	\begin{itemize}
		\item $A_{1,r}\left(\alpha_1,\frac{\alpha_1(1-2\alpha_1)}{2\beta_1},-\frac{\alpha^2_1(1-2\alpha_1)}{2\beta^2_1}-\frac{\alpha_1}{\beta_1},\beta_1\right),$ where $\alpha_1\beta_1\neq 0,$
		\item $A_{1,r}(0,\alpha_2,-2\alpha_2(1+\alpha_2),0),$ where $\alpha_2(1+\alpha_2)\neq 0,$
		\item $A_{1,r}\left(\frac{1}{2},-1,\alpha_4,0\right),\ \ A_{2,r}\left(\frac{1}{2},0,\beta_2\right),\ \ A_{3,r}\left(\frac{1}{2},0,\beta_2\right),\ \ A_{5,r}\left(\frac{1}{2},\beta_2\right).$
	\end{itemize}
\end{theorem}
The results above are represented in the following table, where the units as well are shown.
\noindent\makebox[\linewidth]{
	\begin{tabular}{ |c|c| }
		\hline
		Algebra & $\mathbf{1_\textit{L}} $ \\
		\hline
		$A_{1,r}\left(\alpha_1,\frac{\alpha_1(1-\alpha_1)}{\beta_1}-\frac{1}{2},\frac{\alpha_1(1-\alpha_1)^2}{\beta^2_1}-\frac{1-\alpha_1}{2\beta_1},\beta_1\right)$, where $\beta_1\neq 0$&$
		\frac{-2(1-\alpha_1)}{\beta_1}e_1+2e_2
		$
		\\
		\hline
		$A_{1,r}\left(1,\alpha_2,\frac{\alpha_2(2\alpha_2+1)}{2},0\right)$ &$
		-(1+2\alpha_2)e_1+2e_2
		$
		\\
		\hline
		$A_{2,r}(\alpha_1,0,\alpha_1),$ where $\alpha_1\neq 0$& $
		\frac{1}{\alpha_1}e_1
		$
		\\
		\hline
		$A_{3,r}(\alpha_1,0,\alpha_1),$ where $\alpha_1\neq 0$& $
		\frac{1}{\alpha_1}e_1
		$
		\\
		\hline
		$A_{5,r}(1, 1)$&$e_1+te_2
		,$ where $t \in \mathbb{R}$
		\\
		\hline
		$A_{5,r}(\alpha_1, \alpha_1)$, where $\alpha_1\neq 0,1.$
		& $		\frac{1}{\alpha_1}e_1
		$
		\\
		\hline
		$A_{7,r}\left(\frac{1}{2},0\right)$&$2e_1
		$
		\\
		\hline
		$A_{8,r}\left(\frac{1}{2},0\right)$&$2e_1
		$
		\\
		\hline
		$A_{10,r}(\frac{1}{2})$& $2e_1
		$
		\\
		\hline
		Algebra &$\mathbf{1_\textit{R}} $ \\
		\hline
		$A_{1,r}\left(\alpha_1,\frac{\alpha_1(1-2\alpha_1)}{2\beta_1},\frac{-\alpha^2_1(1-2\alpha_1)}{2\beta^2_1}-\frac{\alpha_1}{\beta_1},\beta_1\right)$, where $\alpha_1\beta_1\neq 0$&$2e_1+	\frac{2\beta_1}{\alpha_1}e_2
		$
		\\
		\hline
		$A_{1,r}(0,\alpha_2, -2\alpha_2(\alpha_2+1),0)$, where $\alpha_2(1+\alpha_2)\neq 0$ &$2e_1+\frac{1}{\alpha_2}e_2
		$
		\\
		\hline
		$A_{1,r}\left(\frac{1}{2},-1,\alpha_4,0\right)$ &  $2e_1
		$
		\\
		\hline
		$A_{2,r}(\frac{1}{2}, 0,\beta_2)$ &$2e_1
		$
		\\
		\hline
		$A_{3,r}(\frac{1}{2}, 0,\beta_2)$ &$2e_1
		$
		\\
		\hline
		$A_{5,r}(\frac{1}{2},0)$&$2e_1+	te_2
		,$ where $t \in \mathbb{R}$
		\\
		\hline
		$A_{5,r}(\frac{1}{2},\beta_2),$ where $\beta_2\neq0$
		&$2e_1
		$
		\\
		\hline
	\end{tabular}
}

\begin{corollary} Up to isomorphism there are only the following nontrivial $2$-dimensional real unital algebras.
	\[A_{2,r}\left(\frac{1}{2},0,\frac{1}{2}\right)=\left(
	\begin{array}{cccc}
	\frac{1}{2} & 0 &0 & 1 \\
	0 & \frac{1}{2} & \frac{1}{2} & 0
	\end{array}\right),\ A_{3,r}\left(\frac{1}{2},0,\frac{1}{2}\right)=\left(
	\begin{array}{cccc}
	\frac{1}{2} & 0 &0 & -1 \\
	0 & \frac{1}{2} & \frac{1}{2} & 0
	\end{array}\right),\]
	\[ A_{5,r}\left(\frac{1}{2},\frac{1}{2}\right)=\left(
	\begin{array}{cccc}
	\frac{1}{2} & 0 &0 & 0 \\
	0 & \frac{1}{2} & \frac{1}{2} & 0
	\end{array}\right).\]\end{corollary}
Among these algebras only $A_{3,r}\left(\frac{1}{2},0,\frac{1}{2}\right)$ is a division algebra and it is isomorphic to the algebra of complex numbers.

\vskip 0.4 true cm

\end{document}